\documentclass[a4paper,12pt]{article}
\usepackage{amssymb}
\usepackage{amsmath,amsfonts,amssymb}
\usepackage{mathrsfs}

\DeclareMathOperator{\trace}{Tr}
\DeclareMathOperator{\Op}{op}
\def\D{\mathbb D}\def\C{\mathbb C}\def\Z{\mathbb Z}\def\R{\mathbb R}\def\T{\mathbb T}
\def\rhodiff{\mbox{$\triangle\!\!\!\!\ast\,$}}
\def\d{{\rm d} }
\def\e{{\rm e} }

\title{On multipliers on compact Lie groups
\thanks{The research was supported by the {\rm EPSRC} grant {\rm EP/G007233/1}.}}
\author{M. V. Ruzhansky and J. Wirth}
\date{}
\begin{document}

 \setlength{\baselineskip}{25pt}

\newtheorem{theorem}{Theorem}
\newtheorem{corollary}[theorem]{Corollary}
\hyphenation{SPbGU}

  \newcommand{\Rn}{{{\mathbb R}^n}}
  \newcommand{\Rone}{{\mathbb R}}

\maketitle

\begin{abstract}\footnote{Keywords: mutipliers, pseudo-differential operators, Lie groups}
In this note we announce $L^p$ multiplier theorems for
invariant and non-invariant operators on compact Lie groups
in the spirit of the well-known H\"ormander-Mikhlin theorem 
on $\mathbb R^n$ and its variants on tori $\mathbb T^n$. 
Applications are given to the mapping properties of 
pseudo-differential operators on $L^p$-spaces 
and to a-priori estimates 
for non-hypoelliptic operators.
\end{abstract}

{\bf 1. Introduction.}

Let $G$ be a compact Lie group of dimension $n$, 
with identity $1$ and the 
unitary dual $\widehat G$.
The following considerations are based on the group Fourier transform
\begin{equation}\label{eq:FI}
{\mathscr F}\phi = \widehat \phi(\xi) = 
\int_G \phi(x) \xi(x)^* \d x,\qquad
  \phi(x) = \sum\nolimits_{[\xi]\in\widehat G} d_\xi 
  \trace(\xi(x) \widehat \phi(\xi) ) = 
  {\mathscr F}^{-1} [\widehat\phi]  
\end{equation}
defined in terms of equivalence classes 
$[\xi]$ of irreducible unitary 
representations $\xi : G \to \mathrm U(d_\xi)$ 
of dimension $d_\xi$. 
The Peter--Weyl theorem on $G$ implies in 
particular that this pair of 
transforms is inverse to each other and that the Plancherel identity
\begin{equation}\label{eq:Plancherel}
   \|\phi\|_2^2 = \sum\nolimits_{[\xi]\in \widehat G} d_\xi \|\widehat\phi(\xi)\|_{HS}^2
   =: \|\widehat{\phi}\|_{\ell^2(\widehat G)}
\end{equation}
holds true for all $\phi\in L^2(G)$. Here $\|\widehat
\phi(\xi)\|_{HS}^2 = \trace(\widehat \phi(\xi)\widehat\phi(\xi)^*)$
denotes the Hilbert--Schmidt norm of matrices. 
The Fourier inversion 
statement \eqref{eq:FI} is valid for all 
$\phi\in\mathcal D'(G)$ and
the Fourier series converges in $C^\infty(G)$ provided $\phi$ is
smooth. It is further convenient to denote
$
  \langle\xi\rangle = \max\{1, \lambda_\xi\},
$
where $-\lambda_\xi^2$ is the eigenvalue of the 
Laplace-Beltrami (Casimir) operator acting on the matrix coefficients 
associated to the representation $\xi$. The Sobolev spaces 
can be characterised by Fourier coefficients as
\begin{equation*}
 \phi\in H^s(G) \quad\Longleftrightarrow \quad 
 \langle\xi\rangle^s \widehat\phi(\xi) \in \ell^2(\widehat G), 
\end{equation*}
where $\ell^2(\widehat G)$ is defined as the space of matrix-valued
sequences such that the sum on the right-hand side of 
\eqref{eq:Plancherel} is finite.

In the following we consider continuous 
linear operators $A:C^\infty(G) \to \mathcal D'(G)$, which can be characterised
by their symbol
\begin{equation}\label{EQ:quantisation0}
  \sigma_A(x,\xi)  = \xi(x)^* (A\xi)(x)
\end{equation}
which is a function on $G\times\widehat G$ taking matrices from $\C^{d_\xi\times d_\xi}$ as
values. As a consequence of \eqref{eq:FI} we obtain that for any given $\phi\in C^\infty(G)$
the distribution $A\phi\in\mathcal D'(G)$ satisfies
\begin{equation}\label{EQ:quantisation}
  A\phi(x) = \sum\nolimits_{[\xi]\in\widehat G} 
  d_\xi \trace(\xi(x) \sigma_A(x,\xi) \widehat \phi(\xi) ).
\end{equation}
We denote the operator $A$ defined by a symbol $\sigma_A$ as $\Op(\sigma_A)$.
This quantisation and its properties have been consistently
developed in \cite{RTbook} and we refer to it for details.
We speak of a {\em Fourier multiplier} if the symbol $\sigma_A(x,\xi)$ is independent of the first
argument. This is equivalent to requiring that $A$ commutes with left translations. It is evident from the Plancherel identity that such an operator is $L^2$-bounded if and only if
$\sup_{[\xi]\in\widehat G} \|\sigma_A(\xi)\|_{op}<\infty$, where $\|\cdot\|_{op}$ denotes the operator norm
on the inner-product space $\C^{d_\xi}$.

In the book \cite{RTbook}, as well as in the paper \cite{RTW10} the authors gave a characterisation 
of H\"ormander type pseudo-differential operators on $G$ in terms of their matrix-valued symbols. 
The symbol classes, as well as the multiplier theorems given below, depend on the so-called difference operators acting on moderate sequences of matrices, i.e., on elements of 
\begin{equation*}
   \Sigma(\widehat G) = \{ \sigma : \xi \mapsto \sigma(\xi) 
   \in \mathbb C^{d_\xi\times d_\xi} : 
   \| \sigma(\xi) \|_{op} \lesssim \langle 
   \xi\rangle^N \textrm { for some $N$}\}.
\end{equation*}
A difference operator $Q$ of order $\ell$ is defined in terms of  a corresponding function
$q\in C^\infty(G)$, which 
vanishes to (at least) $\ell$th order in the identity 
element $1\in G$ via
\begin{equation}\label{EQ:def-dif}
   Q\sigma  = \mathscr F \left(q(x) {\mathscr F}^{-1}\sigma\right).
\end{equation}
Note, that $\sigma\in\Sigma(\widehat G)$ implies 
${\mathscr F}^{-1}\sigma\in\mathcal D'(G)$ and therefore the 
multiplication with a smooth function is well-defined. 
The main idea of introducing such operators is that
applying differences to symbols of Calderon--Zygmund
operators brings an improvement in the behaviour
of $\Op (Q\sigma)$ since we multiply the integral kernel of
$\Op (\sigma)$ by a function vanishing on its singular
set. 

Different collections of difference operators
have been explored in \cite{RTW10} in the 
pseudo-differential setting.
Difference operators of particular interest arise from 
matrix-coefficients of representations. 
For a fixed irreducible representation $\xi_0$ we 
define the (matrix-valued) difference operator
${}_{\xi_0}\mathbb D=
({}_{\xi_0}\mathbb D_{ij})_{i,j=1,\ldots,d_{\xi_0}}$ 
corresponding to the matrix elements of the matrix-valued function
$\xi_0(x)-\mathrm I$, with 
$q_{ij}(x)=\xi_0(x)_{ij}-\delta_{ij}$ in \eqref{EQ:def-dif},
$\delta_{ij}$ the Kronecker delta.
If the representation is fixed, 
we omit the index $\xi_0$. For a sequence of 
difference operators of this type,
$\D_1={}_{\xi_1}\D_{i_1 j_1},
\D_2={}_{\xi_2}\D_{i_2 j_2}, \ldots,
\D_k={}_{\xi_k}\D_{i_k j_k}$, 
with $[\xi_m]\in \widehat G$, $1\leq i_m,j_m\leq d_{\xi_m}$,
$1\leq m\leq k$,
we define
$\D^\alpha=\D_1^{\alpha_1}\cdots \D_k^{\alpha_k}$. In the sequel we will work with 
a collection $\Delta_0$ of representations chosen as follows.
Let $\widetilde{\Delta_{0}}$ be the collection of the irreducible components of the adjoint
representation, so that
${\rm Ad}=(\dim Z(G)) 1\oplus\bigoplus_{\xi\in \widetilde{\Delta_{0}}}\xi$, where $\xi$
are irreducible representations and $1$ is the trivial one-dimensional representation.
In the case when the centre $Z(G)$ of the group is nontrivial, we extend the collection
$\widetilde{\Delta_{0}}$ to some collection $\Delta_{0}$ by adding to $\widetilde{\Delta_{0}}$
a family of irreducible representations such that their direct sum is nontrivial on
$Z(G)$, and such that the function 
\begin{equation*}\label{eq:rhoDef}
\rho^2(x) = \sum\nolimits_{[\xi]\in\Delta_0} \left(
d_\xi-\trace \xi(x)\right) \ge 0
\end{equation*} 
(which vanishes only in $x=1$) would define the square of some distance function on $G$ 
near the identity element. 
Such an extension is always possible, and we denote by $\Delta_{0}$ any such
extension; in the case of the trivial centre we do not have to take
an extension and we set $\Delta_{0}=\widetilde{\Delta_{0}}$.
We denote further by $\rhodiff$ the second order difference operator associated to $\rho^2(x)$,
$\rhodiff={\mathscr F} \rho^2(x) {\mathscr F}^{-1}$.
In the sequel, when we write ${\mathbb D}^\alpha$, we can always
assume that it is composed only of ${}_{\xi_m} {\mathbb D}_{i_m j_m}$ with
$[\xi_m]\in\Delta_0$.

\bigskip
{\bf 2. Main results.}

The following condition \eqref{eq:HM-cond} is a natural
relaxation from the $L^p$-boundedness of
zero order pseudo-differential operators to a 
multiplier theorem and generalises the H\"ormander--Mikhlin 
(\cite{Mih1, Mih2}, \cite{Horm})
theorem
to arbitrary groups.

\begin{theorem}\label{thm:main1} 
Denote by $\varkappa$ be the smallest even 
integer 
larger than $\frac 12\dim G$. 
Let $A: C^\infty(G) \to \mathcal D'(G)$ be left-invariant.  
Assume that its symbol $\sigma_A$ 
satisfies
\begin{equation}\label{eq:HM-cond}
   \|\rhodiff^{\varkappa/2} \sigma_A(\xi)\|_{op} \le C \langle \xi\rangle^{-\varkappa}
        \qquad\text{and}\qquad  
   \|{\mathbb D}^{\alpha} \sigma_A(\xi) \|_{op} 
   \le C_\alpha \langle\xi\rangle^{-|\alpha|}
\end{equation}
for all multi-indices $\alpha$ with $|\alpha|\le \varkappa-1$,
and for all $[\xi]\in \widehat G$. 
 Then the operator $A$ is 
of weak\footnote{The operator $A$ is said to be of weak type $(1,1)$ if there exists a constant
$C>0$ such that for all $\lambda>0$ and $u\in L^{1}(G)$ the inequality
$\mu\{x\in G:|Au(x)|>\lambda\}\leq C\|u\|_{L^{1}(G)}/\lambda$
holds true, where $\mu$ is the Haar measure on $G$.} 
type $(1,1)$ and $L^p$-bounded for all $1<p<\infty$.
\end{theorem}

We now give some particular applications of Theorem~\ref{thm:main1}. The selection is not complete and indicates a few applications which could be derived from the main result. Full proofs can be found in \cite{RW11}.

\begin{theorem}\label{THM:cor1}
Assume that $\sigma_A\in \mathscr S^0_{\rho}(G)$, 
i.e., by definition, it satisfies inequalities
\begin{equation*} 
    \|\mathbb D^\alpha \sigma_A(\xi)\|_{op} 
    \le C_\alpha \langle\xi\rangle^{-\rho|\alpha|},
\end{equation*}
for some $\rho\in[0,1]$ and all $\alpha$. 
Then $A$ defines a bounded operator from\footnote{
Here $W^{p,r}(G)$ stands for the Sobolev space consisting of all distributions $f$ such that
$(I-\mathscr L)^{r/2}f\in L^{p}(G)$, where $\mathscr L$ is a Laplacian 
(Laplace-Betrami operator, Casimir element) on $G$.}
 $W^{p,r}(G)$ to $L^p(G)$ for 
 $r= \varkappa (1-\rho) | \frac 1p-\frac12|$, $\varkappa$ 
 as in Theorem~\ref{thm:main1} and $1<p<\infty$. 
\end{theorem}

The previous statement applies in particular to the 
parametrices constructed in \cite{RTW10}. We will give two examples on the group
$\mathrm{SU}(2)\cong\mathbb S^3$. Let $\mathrm D_1$,
$\mathrm D_2$ and $\mathrm D_3$ be an orthonormal basis of $\mathfrak{su}(2)$. 
Then both, the
sub-Laplacian $\mathcal L_s = \mathrm D_1^2+\mathrm D_2^2$
as well as the 'heat' operator $\mathcal H=\mathrm D_3-\mathrm D_1^2-\mathrm D_2^2$
have a parametrix from\footnote{
The class $\mathrm{op}\mathscr S^{-1}_{\frac12}({\mathbb S^3})$ is defined as the class 
of operators with symbols
$\sigma_{A}$ satisfying the inequalities
$\|\mathbb D^{\alpha} \sigma_{A}(\xi)\|_{op}\leq C_{\alpha}\langle{\xi}\rangle^{-1-|\alpha|/2}.$}
$\mathrm{op}\mathscr S^{-1}_{\frac12}({\mathbb S^3})$
and therefore the sub-elliptic estimates
\begin{equation}\label{eq:sub-ell-Ls}
\|u\|_{W^{p, 1-|\frac1p-\frac12|}(\mathbb S^3)}\leq 
C_p\|\mathcal L_s u\|_{L^p(\mathbb S^3)}
\qquad\text{and}\qquad
\|u\|_{W^{p, 1-|\frac1p-\frac12|}(\mathbb S^3)}\leq 
C_p\|\mathcal H u\|_{L^p(\mathbb S^3)}
\end{equation}
are valid for all $1<p<\infty$.
The following statement concerns operators which 
are neither locally invertible nor locally hypoelliptic. 

\begin{corollary}\label{COR:vfs}
Let $X$ be a left-invariant
real vector field on $G$. 
Then there exists a discrete exceptional set 
$\mathscr C\subset\mathrm i\mathbb R$, 
such that for any complex number 
$c\not\in\mathscr C$ the operator
$X+c$ is invertible with inverse in $\mathrm{op}\mathscr S^{0}_{0}(G)$.
Consequently, the inequality
$$
  \|f\|_{L^p(G)} \le C_p \|(X+c)f\|_{W^{p,\varkappa|\frac1p-\frac12|}(G)}
$$
holds true for all $1<p<\infty$ and all functions $f$ from that 
Sobolev space, with $\varkappa$ as above. 
\end{corollary}

For the particular case $G=\mathrm{SU}(2)$, the exceptional set coincides with the spectrum
of the skew-selfadjoint realisation of $X$ suitably normalised with
respect to the Killing norm, e.g., $\mathscr C = \mathrm i \frac12 \Z$ if $X=\mathrm D_3$.

The H\"ormander multiplier theorem \cite{Horm},
although formulated in $\R^n$, has a natural analogue on the
torus $\T^n$. The assumtions in Theorem
\ref{thm:main1} on the top order difference brings a refinement
of the toroidal multiplier theorem, at least for some dimensions.
If $G=\T^n=\R^n/\Z^n$, the set $\Delta_0$
can be chosen to consist of $2n$ functions $\e^{\pm 2\pi\mathrm i x_j}$, $1\leq j\leq n$.
Consequently, we have that 
$\rho^2(x)=2n-\sum_{j=1}^n 
\left(\e^{2\pi\mathrm i x_j}+\e^{-2\pi\mathrm i x_j}\right)$ 
in \eqref{eq:rhoDef}, and hence
$\rhodiff \sigma(\xi)=2n \sigma(\xi)-\sum_{j=1}^n 
\left(\sigma(\xi+e_j)+\sigma(x-e_j)\right)$,
where $\xi\in\Z^n$ and $e_j$ is its $j$th unit basis vector
in $\Z^n$.

A (translation) invariant operator $A$ and its symbol $\sigma_A$
are related
 by $\sigma_A(k)=\e^{-2\pi\mathrm i x\cdot k}
(A\e^{2\pi\mathrm i x\cdot k})=(A\e^{2\pi\mathrm i x\cdot k})|_{x=0}$ and
$A\phi(x)=\sum_{k\in\Z^n} e^{2\pi\mathrm i x\cdot k} 
\sigma_A(k) \widehat{\phi}(k)$.
Thus, it follows from Theorem \ref{thm:main1} that, for example
on ${\mathbb T}^3$, a translation invariant operator
$A$ is bounded on $L^p({\mathbb T}^3)$ provided 
that there is a constant
$C>0$ such that $|\sigma_A(k)|\leq C$,
$|k||\sigma_A(k+e_j)-\sigma_A(k)|\leq C$ and
\begin{equation}\label{EQ:T3}
|k|^2| \sigma_A(k)-{\textstyle\frac16}\sum\nolimits_{j=1}^3 
\left(\sigma_A(k+e_j)+\sigma_A(k-e_j)\right)|\leq C, 
\end{equation} 
for all
$k\in\Z^3$ and all (three) unit vectors $e_j$, $j=1,2,3$.
Here we do not make assumptions on all second order
differences in \eqref{EQ:T3}, but only on one of them.

Finally, Theorem~\ref{thm:main1} also implies a boundedness 
statement for operators of form \eqref{EQ:quantisation}. 
Let for this 
$\partial_{x_j}$, $1\le j\le n$, be a 
collection of left invariant first order differential
operators corresponding to some linearly independent family
of the left-invariant vector fields on $G$. As usual, we denote
$\partial_x^{\beta}=\partial_{x_1}^{\beta_1}\cdots
\partial_{x_n}^{\beta_n}$.

\begin{theorem}\label{thm:main3} 
Denote by $\varkappa$ be the smallest even 
integer  
larger than $\frac n2$, 
$n$ the dimension of the group $G$. 
Let $1<p<\infty$ and let $l>\frac{n}{p}$ be an integer.
Let $A: C^\infty(G) \to \mathcal D'(G)$ 
be a linear continuous operator
such that its matrix symbol $\sigma_A$ 
satisfies
\begin{equation*}\label{eq:HM-cond-p}
   \|\partial_x^\beta{\mathbb D}^{\alpha} 
   \sigma_A(x,\xi) \|_{op} \le 
   C_{\alpha,\beta} \langle\xi\rangle^{-|\alpha|}
\end{equation*}
for all multi-indices $\alpha,\beta$ with 
$|\alpha|\le \varkappa$ and 
$|\beta|\leq l$, for all $x\in G$ and $[\xi]\in \widehat G$.  
Then the operator $A$ is 
bounded on $L^p(G)$.
\end{theorem}

{\bf 3. Discussion.} 1.
The conditions are needed for the week type $(1,1)$ property. 
Interpolation allows to reduce assumptions on the number of 
differences for $L^p$-boundedness. The result generalises the corresponding statements 
in the case of the group SU(2) in
\cite{CdG70}, \cite{CW70}, also presented in \cite{CWbook}.

2. Examples similar to \eqref{eq:sub-ell-Ls}
can be given for arbitrary compact 
Lie groups. The assumptions of Theorem \ref{THM:cor1} concerning
the numbers of difference operators can be relaxed to the same
as those in Theorem \ref{thm:main1}.

3. If the operator $A\in \Psi^0(G)$ is the usual pseudo-differential operator
of H\"ormander type of order $0$ on $G$
(i.e. in all local coordinate it belongs to H\"ormander class
$\Psi^{0}(\Rn)$), it was shown in \cite{RTbook} the
estimates \eqref{eq:HM-cond-p} hold for all 
$\alpha,\beta$. The converse is also true. Namely, if 
estimates \eqref{eq:HM-cond-p} hold for all 
$\alpha,\beta$, then $A\in \Psi^0(G)$, cf. \cite{RTW10}.

4. Noncommutative matrix quantisation \eqref{EQ:quantisation0}--\eqref{EQ:quantisation}
has a full symbolic calculus (compositions, adjoints, parametrix, etc.),
which have been established in the monograph \cite{RTbook}.

5. On SU(2) the operators corresponding to our difference operators but defined
explicitly in terms of the Clebsch-Gordan coefficients have been used in 
\cite{CW70,CWbook}. The general definition \eqref{EQ:def-dif}, the main tool
in the present investigation, has been introduced and analysed in
\cite{RTbook} and \cite{RTW10}.


\medskip

 \noindent
{\small
M. V. Ruzhansky, {\rm Imperial College London}, {\rm
ruzh@ic.ac.uk}  \\
J. Wirth, University of Stuttgart, 
{\rm jens.wirth@mathematik.uni-stuttgart.de}}

\end{document}